# EQUIVALENT AND ABSOLUTELY CONTINUOUS MEASURE CHANGES FOR JUMP-DIFFUSION PROCESSES

By Patrick Cheridito[1], Damir Filipović and Marc Yor

*Princeton University, University of Munich and Université Pierre et Marie Curie, Paris VI*

We provide explicit sufficient conditions for absolute continuity and equivalence between the distributions of two jump-diffusion processes that can explode and be killed by a potential.

**1. Introduction.** The purpose of this paper is to give explicit, easy-to-check sufficient conditions for the distributions of two jump-diffusion processes to be equivalent or absolutely continuous. We consider jump-diffusions that can explode and be killed by a potential. These processes are, in general, not semimartingales. We characterize them by their infinitesimal generators.

The structure of the paper is as follows. In Section 2 we introduce the notation and state the paper's main result, which gives sufficient conditions for the distributions of two jump-diffusions to be equivalent or absolutely continuous. The conditions consist of local bounds on the transformation of one generator into the other one and the assumption that the martingale problem for the second generator has for all initial distributions a unique solution. The formulation of the main theorem involves two sequences of stopping times. Stopping times of the first sequence stop the process before it explodes. The second sequence consists of exit times of the process from regions in the state space where the transformation of the first generator into the second one can be controlled. Our main result applies also in situations where the generalized Novikov condition ([19], Théorème IV.3) or Kazamaki-like criteria (e.g., [14, 15, 16]) are not satisfied. In Section 3 we show how $X$ can be turned into a semimartingale by embedding it in a larger state space and stopping it before it explodes. The results of Section 3 are needed in the

Received March 2004; revised October 2004.
[1]Supported by the Swiss National Science Foundation and Credit Suisse.

proof of the paper's main theorem, which is given in Section 4. In Section 5 we prove a stronger version of the result of Section 2 for a particular set-up, involving the carré-du-champ operator. In Section 6 we illustrate the main result by showing how the characteristics of a Cox–Ingersoll–Ross [3] short rate process with additional jumps and a potential can be altered by an absolutely continuous or equivalent change of measure.

There exists a vast literature on the absolute continuity of stochastic processes, and below we quote some related publications. In contrast to many of those works, the primary goal of this paper is to provide results that are based on explicit assumptions which are easy to verify in typical applications. For two applications in finance, see [2] and [1], which contain measure changes for multi-dimensional diffusion models and multi-dimensional jump-diffusion models with explosion and potential, respectively.

Itô and Watanabe [10], Kunita [18] and Palmowski and Rolski [24] discuss absolute continuity for general classes of Markov processes.

Kunita [17] characterizes the class of all absolutely continuous Markov processes with respect to a given Markov process. A special discussion for Lévy processes can be found in Sato [29], Section 33.

Dawson [4], Liptser and Shiryaev [21] Kabanov, Liptser and Shiryaev [12], Rydberg [28] and Hobson and Rogers [9] discuss absolute continuity of solutions to stochastic differential equations. They are similar in spirit to Kadota and Shepp [13], which contains sufficient conditions for the distribution of a Brownian motion with stochastic drift to be absolutely continuous with respect to the Wiener measure.

Pitman and Yor [25] and Yor [33] study mutual absolute continuity of (squared) Bessel processes.

Lepingle and Mémin [19] and Kallsen and Shiryaev [14] provide conditions for the uniform integrability of exponential local martingales in a general semimartingale framework (see also Remark 2.7 below), extending the classical results by Novikov [23] and Kazamaki [15].

Discussions of measure changes in a finance context can be found in Sin [30], Lewis [20], Delbaen and Shirakawa [5, 6].

Wong and Heyde [32] give necessary and sufficient conditions for the stochastic exponential of a Brownian motion integral to be a martingale in terms of the explosion time of an associated process.

Among various excellent text books that discuss changes of measure in varying degree of generality are, for example, McKean [22], Rogers and Williams [26], Jacod and Shiryaev [11], Revuz and Yor [27] and Strook [31].

**2. Statement of the main result.** Let $E$ be a closed subset of $\mathbb{R}^d$ and $E_\Delta = E \cup \{\Delta\}$ the one-point compactification of $E$. If not mentioned otherwise, any measurable function $f$ on $E$ is extended to $E_\Delta$ by setting $f(\Delta) := 0$. We let $\Omega$ be the space of càdlàg functions $\omega : \mathbb{R}_+ \to E_\Delta$ such that $\omega(t-) = \Delta$ or $\omega(t) = \Delta$ implies $\omega(s) = \Delta$ for all $s \geq t$. $(X_t)_{t \geq 0}$ is the coordinate process, given by

$$X_t(\omega) := \omega(t), \qquad t \geq 0.$$

It generates the $\sigma$-algebra,

$$\mathcal{F}^X := \sigma(X_s : s \geq 0),$$

and the filtration

$$\mathcal{F}^X_t := \sigma(X_s : 0 \leq s \leq t), \qquad t \geq 0.$$

It follows from Proposition 2.1.5 (a) in [8] that, for all closed subsets $\Gamma$ of $E_\Delta$,

$$\inf\{t \mid X_{t-} \in \Gamma \text{ or } X_t \in \Gamma\} \qquad \text{is an } (\mathcal{F}^X_t)\text{-stopping time.}$$

Hence,

$$T_\Delta := \inf\{t \mid X_t = \Delta\} = \inf\{t \mid X_{t-} = \Delta \text{ or } X_t = \Delta\}$$

is an $(\mathcal{F}^X_t)$-stopping time. Note that

$$X_\cdot = \Delta \qquad \text{on } [T_\Delta, \infty)$$

so that $T_\Delta$ can be viewed as the lifetime of $X$. For the handling of explosion, we introduce the $(\mathcal{F}^X_t)$-stopping times

$$T'_n := \inf\{t \mid \|X_{t-}\| \geq n \text{ or } \|X_t\| \geq n\}, \qquad n \geq 1,$$

where $\|\cdot\|$ denotes the Euclidean norm on $\mathbb{R}^d$ and $\|\Delta\| := \infty$. Clearly, $T'_n \leq T_\Delta$, for all $n \geq 1$. A transition to $\Delta$ occurs either by a jump or by explosion. Accordingly, we define the $(\mathcal{F}^X_t)$-stopping times

$$T_{\text{jump}} := \begin{cases} T_\Delta, & \text{if } T'_n = T_\Delta \text{ for some } n, \\ \infty, & \text{if } T'_n < T_\Delta \text{ for all } n, \end{cases}$$

$$T_{\text{expl}} := \begin{cases} T_\Delta, & \text{if } T'_n < T_\Delta \text{ for all } n, \\ \infty, & \text{if } T'_n = T_\Delta \text{ for some } n, \end{cases}$$

$$T_n := \begin{cases} T'_n, & \text{if } T'_n < T_\Delta, \\ \infty, & \text{if } T'_n = T_\Delta. \end{cases}$$

Note that $\{T_{\text{jump}} < \infty\} \cap \{T_{\text{expl}} < \infty\} = \varnothing$, $\lim_{n \to \infty} T_n = T_{\text{expl}}$, and $T_n < T_{\text{expl}}$ on $\{T_{\text{expl}} < \infty\}$. Hence, $T_{\text{expl}}$ is predictable with announcing sequence $T_n \wedge n$ (see [11], I.2.15.a).

Since, by definition, $\Omega$ contains only paths that stay in $\Delta$ after explosion or after a jump to $\Delta$, the filtration $(\mathcal{F}^X_t)$ has the property stated in Proposition 2.1 below, whose proof is given in the Appendix.

PROPOSITION 2.1.  *Let $T$ be an arbitrary $(\mathcal{F}_t^X)$-stopping time. Then*

$$\mathcal{F}_T^X = \mathcal{F}_{T \wedge T_{\mathrm{expl}}}^X = \sigma\left(\bigcup_{n \geq 1} \mathcal{F}_{T \wedge T_n}^X\right).$$

Fix a bounded and continuous function $\chi : \mathbb{R}^d \to \mathbb{R}^d$ such that $\chi(\xi) = \xi$ on a neighborhood of 0. Let $\alpha, \beta, \gamma$ be measurable mappings on $E$ with values in the set of positive semi-definite symmetric $d \times d$-matrices, $\mathbb{R}^d$ and $\mathbb{R}_+$, respectively. Furthermore, let $\mu$ be a transition kernel from $E$ to $\mathbb{R}^d$ and assume that the functions

(2.1)
$$\alpha(\cdot), \beta(\cdot), \gamma(\cdot) \text{ and } \int_{\mathbb{R}^d} (\|\xi\|^2 \wedge 1) \mu(\cdot, d\xi)$$

are bounded on every compact subset of $E$.

Then,

$$\mathcal{A}f(x) := \frac{1}{2} \sum_{i,j=1}^d \alpha_{ij}(x) \frac{\partial^2 f(x)}{\partial x_i \, \partial x_j} + \sum_{i=1}^d \beta_i(x) \frac{\partial f(x)}{\partial x_i} - \gamma(x) f(x)$$
$$+ \int_{\mathbb{R}^d} (f(x + \xi) - f(x) - \langle \nabla f(x), \chi(\xi) \rangle) \mu(x, d\xi)$$

defines a linear operator from the space of $C^2$-functions on $E$ with compact support, $C_c^2(E)$, to the space of bounded measurable functions on $E$, $B(E)$.

DEFINITION 2.2.  We say that a probability measure $\mathbb{P}$ on $(\Omega, \mathcal{F}^X)$ is a solution of the martingale problem for $\mathcal{A}$ if, for all $f \in C_c^2(E)$,

$$M_t^f := f(X_t) - f(X_0) - \int_0^t \mathcal{A}f(X_s) \, ds, \qquad t \geq 0,$$

is a $\mathbb{P}$-martingale with respect to $(\mathcal{F}_t^X)$. We say that the martingale problem for $\mathcal{A}$ is well posed if for every probability distribution $\eta$ on $E$, there exists a unique solution $\mathbb{P}$ of the martingale problem for $\mathcal{A}$ such that $\mathbb{P} \circ X_0^{-1} = \eta$.

REMARK 2.3.  1. If $\mathbb{P}$ is a solution of the martingale problem for $\mathcal{A}$, then with respect to $\mathbb{P}$, $X$ is a possibly nonconservative, time-homogenous jump-diffusion process. The time-homogeneous case can be included in the above set-up by identifying one component of $X$ with time $t$.

2. If $\mathbb{P}$ is a solution of the martingale problem for $\mathcal{A}$, then $M^f$ is, for all $f \in C_c^2(E)$, also a $\mathbb{P}$-martingale with respect to $(\mathcal{F}_{t+}^X)$. Indeed, since all paths of $M^f$ are right-continuous, it follows from the backwards martingale convergence theorem that, for all $t, s \in \mathbb{R}$ such that $t < s$,

$$\mathbb{E}_{\mathbb{P}}[M_s^f | \mathcal{F}_{t+}^X] = \lim_{r \searrow t} \mathbb{E}_{\mathbb{P}}[M_s^f | \mathcal{F}_r^X] = \lim_{r \searrow t} M_r^f = M_t^f.$$

3. It is easy to see that if the martingale problem for $\mathcal{A}$ is well posed, then for every probability distribution $\eta$ on $E_\Delta$, there exists a unique solution $\mathbb{P}$ of the martingale problem for $\mathcal{A}$ such that $\mathbb{P} \circ X_0^{-1} = \eta$.

4. Throughout, we make use of the fact that $\int_0^t f(X_{u-})\, dS_u = \int_0^t f(X_u)\, dS_u$, for a continuous semimartingale $S$ and every measurable function $f$ such that the integrals are defined.

Let $\tilde{\mathcal{A}}$ be a second linear operator from $C_c^2(E)$ to $B(E)$, given by

$$\tilde{\mathcal{A}}f(x) := \frac{1}{2} \sum_{i,j=1}^{d} \alpha_{ij}(x) \frac{\partial^2 f(x)}{\partial x_i\, \partial x_j} + \sum_{i=1}^{d} \tilde{\beta}_i(x) \frac{\partial f(x)}{\partial x_i} - \tilde{\gamma}(x) f(x)$$

(2.2)
$$+ \int_{\mathbb{R}^d} (f(x+\xi) - f(x) - \langle \nabla f(x), \chi(\xi) \rangle) \tilde{\mu}(x, d\xi),$$

where $\tilde{\beta}$ and $\tilde{\gamma}$ are measurable mappings from $E$ to $\mathbb{R}^d$ and $\mathbb{R}_+$, respectively, and $\tilde{\mu}$ is a transition kernel from $E$ to $\mathbb{R}^d$ such that $\tilde{\beta}$, $\tilde{\gamma}$ and $\tilde{\mu}$ satisfy the condition (2.1).

Let $U$ be an open subset of $E$, that is, $U = U' \cap E$ for some open subset $U'$ of $\mathbb{R}^d$. Assume that there exist measurable mappings

$$\phi_1 : U \to \mathbb{R}^d, \qquad \phi_2 : U \to (0, \infty) \quad \text{and} \quad \phi_3 : U \times \mathbb{R}^d \to (0, \infty)$$

such that, for all $x \in U$,

(2.3)
$$\tilde{\beta}(x) = \beta(x) + \alpha(x)\phi_1(x) + \int_{\mathbb{R}^d} (\phi_3(x, \xi) - 1) \chi(\xi) \mu(x, d\xi),$$
$$\tilde{\gamma}(x) = \phi_2(x) \gamma(x),$$
$$\tilde{\mu}(x, d\xi) = \phi_3(x, \xi) \mu(x, d\xi).$$

Let $U^1 \subset U^2 \subset \cdots$ be an increasing sequence of open subsets of $E$ such that $U = \bigcup_{n \geq 1} U^n$. We denote $U_\Delta = U \cup \{\Delta\}$ and $U_\Delta^n = U^n \cup \{\Delta\}$, $n \geq 1$. For all $n \geq 1$, we define

$$R_n := \inf\{t \mid X_{t-} \notin U_\Delta^n \text{ or } X_t \notin U_\Delta^n\}.$$

Note that

$$R_n = \begin{cases} R_n', & \text{if } R_n' < T_\Delta, \\ \infty, & \text{if } R_n' = T_\Delta, \end{cases}$$

where

$$R_n' := \inf\{t \mid X_{t-} \notin U^n \text{ or } X_t \notin U^n\}, \qquad n \geq 1.$$

Since the sets $U^n$ are open in the topology of $E_\Delta$, it follows from Proposition 2.1.5(a) of [8] that all $R'_n$, $R_n$ and therefore also,

$$R_\infty := \lim_{n \to \infty} R_n = \inf\{t \mid X_{t-} \notin U_\Delta \text{ or } X_t \notin U_\Delta\},$$
$$S_n := R_n \wedge T_n \wedge n, n \geq 1,$$
$$S_\infty := \lim_{n \to \infty} S_n = R_\infty \wedge T_{\text{expl}}$$

are $(\mathcal{F}_t^X)$-stopping times. While the sequence $T_1 \leq T_2 \leq \cdots$ takes care of a possible explosion of $X$, the sequence $S_1 \leq S_2 \leq \cdots$ appropriately localizes the stochastic logarithm of the density process for the measure change, see (4.2) below. In view of (2.1) and the convention $f(\Delta) = 0$ for measurable functions $f$,

$$\Lambda_n := \tfrac{1}{2} \int_0^{S_n} \langle \alpha(X_s)\phi_1(X_s), \phi_1(X_s) \rangle \, ds$$
$$+ \int_0^{S_n} (\phi_2(X_s) \log \phi_2(X_s) - \phi_2(X_s) + 1) \gamma(X_s) \, ds$$
$$+ \int_0^{S_n} \int_{\mathbb{R}^d} (\phi_3(X_s, \xi) \log \phi_3(X_s, \xi) - \phi_3(X_s, \xi) + 1) \mu(X_s, d\xi) \, ds$$

is well defined for all $n \geq 1$. With this notation we have the following:

THEOREM 2.4. *Let $\mathbb{P}$ be a solution of the martingale problem for $\mathcal{A}$ and $\mathbb{Q}$ a solution of the martingale problem for $\tilde{\mathcal{A}}$ such that $\mathbb{Q}|_{\mathcal{F}_0^X} \ll \mathbb{P}|_{\mathcal{F}_0^X}$. Assume that for $\tilde{\mathcal{A}}$, the martingale problem is well posed and that*

(2.4) $$\mathbb{E}_\mathbb{P}[e^{\Lambda_n}] < \infty,$$

*for all $n \geq 1$.*

*Then there exists a nonnegative càdlàg $\mathbb{P}$-supermartingale $(D_t)_{t \geq 0}$ such that, for any $(\mathcal{F}_t^X)$-stopping time $T$, the following properties hold:*

(2.5)  $\mathbb{Q}|_{\mathcal{F}_T^X \cap \{T < S_\infty\}} = D_T \cdot \mathbb{P}|_{\mathcal{F}_T^X \cap \{T < S_\infty\}}.$

(2.6)  *If $\mathbb{Q}[T < S_\infty] = 1$, then $\mathbb{Q}|_{\mathcal{F}_T^X} = D_T \cdot \mathbb{P}|_{\mathcal{F}_T^X}.$*

(2.7)  *If $\mathbb{Q}|_{\mathcal{F}_0^X} \sim \mathbb{P}|_{\mathcal{F}_0^X}$ and $\mathbb{P}[T < S_\infty] = \mathbb{Q}[T < S_\infty] = 1$,*

*then $\mathbb{Q}|_{\mathcal{F}_T^X} \sim \mathbb{P}|_{\mathcal{F}_T^X}.$*

(2.8)  *If $\mathbb{Q}[T < R_\infty] = 1$ and $(D_{T \wedge T_n})_{n \geq 1}$ is $\mathbb{P}$-uniformly integrable,*

*then $\mathbb{Q}|_{\mathcal{F}_T^X} = D_T \cdot \mathbb{P}|_{\mathcal{F}_T^X}.$*

REMARK 2.5. The following is an easy-to-check sufficient criterion for (2.4): Assume that for every $n \geq 1$, there exists a finite constant $K_n$ such that, for all $x \in U^n$,

$$\langle \alpha(x)\phi_1(x), \phi_1(x)\rangle \leq K_n, \tag{2.9}$$

$$(\phi_2(x)\log\phi_2(x) - \phi_2(x) + 1)\gamma(x) \leq K_n, \tag{2.10}$$

$$\int_{\mathbb{R}^d}(\phi_3(x,\xi)\log\phi_3(x,\xi) - \phi_3(x,\xi) + 1)\mu(x,d\xi) \leq K_n. \tag{2.11}$$

Then (2.4) is satisfied.

REMARK 2.6. If $\mathbb{P}[S_\infty = \infty] = 1$, we obtain from (2.5) the loss of mass of the $\mathbb{P}$-supermartingale $(D_t)_{t \geq 0}$

$$1 - \mathbb{E}_\mathbb{P}[D_t] = 1 - \mathbb{Q}[t < S_\infty] = \mathbb{Q}[S_\infty \leq t], \qquad t \in [0, \infty).$$

REMARK 2.7. For $\phi_3(x,\xi) = e^{\langle\phi_1(x),\xi\rangle}$, our measure changes are of the same form as the generalized Esscher transforms discussed in [14] (see Theorem 2.19 in [14] or Theorem III.7.23 in [11]).

**3. Turning $X$ into a semimartingale.** In this section we show some preliminary results that we will need in the proof of Theorem 2.4. The notation is the same as in Section 2. For any process $Y$ and stopping time $T$, we denote by $Y^T$ the stopped process given by $Y_t^T := Y_{t \wedge T}$, $t \geq 0$.

Assume that $\mathbb{P}$ is a solution of the martingale problem for $\mathcal{A}$. Since the coordinate process can explode and be killed, it is, in general, not a semimartingale with respect to $\mathbb{P}$. To turn it into a semimartingale, we stop it before it explodes and identify the state $\Delta$ with an arbitrary point $\partial$ in $\mathbb{R}^d \setminus E$. Without loss of generality, we can assume that such a point exists. If $E = \mathbb{R}^d$, we embed $E$ in $\mathbb{R}^{d+1}$ by the map $(x_1, \ldots, x_d) \mapsto (x_1, \ldots, x_d, 0)$ and adjust $\alpha$, $\beta$, $\mu$ and $\chi$ as follows: For all $x \in E$, we extend $\alpha(x)$ to a $(d+1) \times (d+1)$-matrix by setting $\alpha(x)_{i,d+1} = \alpha(x)_{d+1,i} := 0$ for all $i = 1, \ldots, d+1$. $\beta(x)$ is elongated to a $(d+1)$-dimensional vector by $\beta(x)_{d+1} := 0$. The measure $\mu(x, \cdot)$ is extended to $\mathbb{R}^{d+1}$ by defining $\mu(x, \mathbb{R}^{d+1} \setminus \mathbb{R}^d) := 0$. Finally, the truncation function $\chi$ can be extended to a bounded and continuous function from $\mathbb{R}^{d+1}$ to $\mathbb{R}^{d+1}$ such that $\chi(\xi) = \xi$ on a neighborhood of $0$, or simply by setting it equal to zero on $\mathbb{R}^{d+1} \setminus \mathbb{R}^d$. Then, a probability measure $\mathbb{P}$ on $(\Omega, \mathcal{F}^X)$ is a solution of the martingale problem for $\mathcal{A}$ in the $\mathbb{R}^{d+1}$-framework if and only if it is in the $\mathbb{R}^d$-framework.

The process

$$\hat{X} := X\mathbb{1}_{[0,T_\Delta)} + \partial\mathbb{1}_{[T_\Delta,\infty)}$$

is also $(\mathcal{F}_t^X)$-adapted and has right-continuous paths in $\mathbb{R}^d$. However, $\hat{X}_{T_\Delta-} = \Delta$ (explosion) is still possible for this process.

Let $T$ be a $(\mathcal{F}_t^X)$-stopping time such that $T < T_{\text{expl}}$, then

$$\bigcup_{n \geq 1} \{T < T_n\} = \Omega, \tag{3.1}$$

and, therefore, (2.1) implies that the following $(\mathcal{F}_t^X)$-predictable processes and random measure are well defined for all $\omega$:

$$B_t^T := \int_0^{t \wedge T} \beta(X_s) + \gamma(X_s)\chi(\partial - X_s)\,ds,$$

$$C_t^T := \int_0^{t \wedge T} \alpha(X_s)\,ds,$$

$$\nu^T(dt, d\xi) := [\mu(X_t, d\xi) + \gamma(X_t)\delta_{\partial - X_t}(d\xi)]\mathbb{1}_{\{t \leq T\}}\,dt.$$

Condition (2.1) also guarantees that $\nu^T$ satisfies Condition 2.13 on page 77 of [11]. Note that one can choose $\chi$ with compact support such that $\chi(\partial - x) = 0$ for all $x \in E$. In that case, the expression for $B^T$ becomes simpler. For $f \in C_b^2(\mathbb{R}^d)$ (the space of bounded $C^2$-functions on $\mathbb{R}^d$), define the process

$$\mathcal{M}^{f,T} := f(\hat{X}^T) - f(\hat{X}_0^T) - \frac{1}{2}\sum_{i,j=1}^d \frac{\partial^2 f(\hat{X}^T)}{\partial x_i \partial x_j} \cdot C_{ij}^T - \nabla f(\hat{X}^T) \cdot B^T$$

$$- (f(\hat{X}^T + \xi) - f(\hat{X}^T) - \langle \nabla f(\hat{X}^T), \chi(\xi)\rangle) * \nu^T$$

("$\cdot$" denotes stochastic integration with respect to a semimartingale and "$*$" stochastic integration with respect to a random measure, for the definition of stochastic integrals with respect to semimartingales and random measures, see, e.g., [11]). The restriction of a function $f \in C_c^2(\mathbb{R}^d)$ to $E$ is in $C_c^2(E)$. Recall that by convention, $f(\Delta) = \alpha(\Delta) = \beta(\Delta) = \gamma(\Delta) = \mu(\Delta, \cdot) = 0$. Thus, it can easily be checked that

$$\mathcal{M}_t^{f,T} = M_t^{f,T} + f(\partial)N_t^T, \qquad t \geq 0, \tag{3.2}$$

where

$$N_t^T := \mathbb{1}_{\{0 < T_\Delta \leq t \wedge T\}} - \int_0^{t \wedge T} \gamma(X_s)\,ds, \qquad t \geq 0.$$

LEMMA 3.1. *Let $T$ be an $(\mathcal{F}_t^X)$-stopping time with $T < T_{\text{expl}}$. Then the process $N^T$ is an $((\mathcal{F}_{t+}^X), \mathbb{P})$-martingale.*

PROOF. Fix $n \geq 1$. We first show that $N^{T_n}$ is an $((\mathcal{F}_{t+}^X), \mathbb{P})$-martingale. Let $(f_k)$ be a sequence in $C_c^2(\mathbb{R}^d)$ with $0 \leq f_k \leq 1$ and $f_k = 1$ on the ball with center 0 and radius $k$, $B_k$. By Remark 2.3 part 2, $M^{f_k,T_n}$ is an $((\mathcal{F}_{t+}^X), \mathbb{P})$-martingale for every $k$. Note that $T_n = 0$ if $\|X_0\| \geq n$. Hence, we have, for all $k > n$,

$$\begin{aligned} M_t^{f_k,T_n} &= f_k(X_t^{T_n}) - f_k(X_0) - \int_0^{t \wedge T_n} \mathcal{A} f_k(X_s) \, ds \\ &= f_k(X_t^{T_n}) - f_k(X_0) \\ &\quad + \int_0^{t \wedge T_n} \left( \gamma(X_s) - \int_{\mathbb{R}^d \setminus B_{k-n}} (f_k(X_s + \xi) - 1) \mu(X_s, d\xi) \right) ds. \end{aligned}$$

Clearly, for all $\omega$,

$$\lim_{k \to \infty} f_k(X_{t \wedge T_n}) = \mathbb{1}_{\{t \wedge T_n < T_\Delta\}}.$$

Moreover, it can be deduced from (2.1) and Lebesgue's dominated convergence theorem that, for all $\omega$,

$$\int_0^{t \wedge T_n} \int_{\mathbb{R}^d \setminus B_{k-n}} |f_k(X_s + \xi) - 1| \mu(X_s, d\xi) \, ds \leq \int_0^{t \wedge T_n} \mu(X_s, \mathbb{R}^d \setminus B_{k-n}) \, ds \to 0,$$

as $k \to \infty$. Furthermore, it follows from (2.1) that there exists a constant $c_n$ such that

$$|M_t^{f_k,T_n}| \leq 1 + \int_0^{t \wedge T_n} (|\gamma(X_s)| + \mu(X_s, \mathbb{R}^d \setminus B_{k-n})) \, ds \leq 1 + c_n t,$$

for all $k \geq n$. Hence, it follows from Lebesgue's dominated convergence theorem that for all $t \geq 0$,

$$-M_t^{f_k,T_n} \to N_t^{T_n} \qquad \text{in } L^1 \text{ as } k \to \infty,$$

which shows that $N^{T_n}$ is an $((\mathcal{F}_{t+}^X), \mathbb{P})$-martingale. This and (3.1) imply that $N^T$ is an $((\mathcal{F}_{t+}^X), \mathbb{P})$-local martingale, and, therefore, by the Doob–Meyer decomposition theorem ([11], I.3.15), $N^T$ is also a uniformly integrable martingale with respect to $((\mathcal{F}_{t+}^X), \mathbb{P})$. □

Notice that $T < T_{\text{expl}}$ implies $\{T_\Delta \leq t \wedge T\} = \{T_{\text{jump}} \leq t \wedge T\}$. Hence, Lemma 3.1 says that $\int_0^{t \wedge T} \gamma(X_s) \, ds$ is the predictable compensator for the time of a jump of the stopped process $X^T$ to $\Delta$. As a consequence, we obtain that $T_{\text{jump}} = \infty$ $\mathbb{P}$-almost surely on $\{X_0 \neq \Delta\}$ if and only if $\gamma(X_t) = 0$ $\mathbb{P}$-almost surely for all $t$.

PROPOSITION 3.2. *Assume that $\mathbb{P}$ is a solution of the martingale problem for $\mathcal{A}$ and $T$ is an $(\mathcal{F}_t^X)$-stopping time such that $T < T_{\text{expl}}$. Then for all $f \in C_b^2(\mathbb{R}^d)$, $\mathcal{M}^{f,T}$ is a local martingale on $(\Omega, (\mathcal{F}_{t+}^X)_{t \geq 0}, \mathbb{P})$ and $\hat{X}^T$ is a semimartingale on $(\Omega, (\mathcal{F}_{t+}^X)_{t \geq 0}, \mathbb{P})$ with characteristics $(B^T, C^T, \nu^T)$ with respect to the truncation function $\chi$.*

PROOF. Fix $n \geq 1$. In view of (3.2), Remark 2.3 part 2 and Lemma 3.1, $\mathcal{M}^{f,T_n}$ is an $((\mathcal{F}_{t+}^X), \mathbb{P})$-martingale for all $f \in C_c^2(\mathbb{R}^d)$.

Now let $f \in C_b^2(\mathbb{R}^d)$. Then $f f_k \in C_c^2(\mathbb{R}^d)$, where the $f_k \in C_c^2(\mathbb{R}^d)$ are as in the proof of Lemma 3.1, and for all $k \geq n$,

$$|\mathcal{M}_t^{f,T_n} - \mathcal{M}_t^{f f_k, T_n}|$$
$$\leq |f(\hat{X}_t^{T_n}) - f f_k(\hat{X}_t^{T_n})|$$
$$+ \int_0^{t \wedge T_n} \int_{\mathbb{R}^d \setminus B_{k-n}} |f(X_s + \xi) - f f_k(X_s + \xi)| \nu^{T_n}(ds, d\xi)$$
$$\leq |f(\hat{X}_{t \wedge T_n}) - f f_k(\hat{X}_{t \wedge T_n})| + \|f\|_\infty \int_0^{t \wedge T_n} \nu^{T_n}(ds, \mathbb{R}^d \setminus B_{k-n}).$$

Obviously,

$$|f(\hat{X}_{t \wedge T_n}) - f f_k(\hat{X}_{t \wedge T_n})| \to 0 \quad \text{in } L^1 \text{ as } k \to \infty,$$

and as in the proof of Lemma 3.1, it can be deduced from (2.1) that

$$\int_0^{t \wedge T_n} \nu^{T_n}(ds, \mathbb{R}^d \setminus B_{k-n}) \to 0 \quad \text{in } L^1 \text{ as } k \to \infty.$$

Hence, $\mathcal{M}^{f,T_n}$ is an $((\mathcal{F}_{t+}^X), \mathbb{P})$-martingale, for all $n \geq 1$. This, together with (3.1), implies that $\mathcal{M}^{f,T}$ is an $((\mathcal{F}_{t+}^X), \mathbb{P})$-local martingale. Thus, it follows from [11], II.2.42, that $\hat{X}^T$ is an $((\mathcal{F}_{t+}^X), \mathbb{P})$-semimartingale with the claimed characteristics.

□

**4. Proof of Theorem 2.4.** There exists a nonnegative, $\mathcal{F}_0^X$-measurable random variable $D_0$ such that

$$\mathbb{Q}|_{\mathcal{F}_0^X} = D_0 \cdot \mathbb{P}|_{\mathcal{F}_0^X}.$$

For each $n \geq 1$, let $\hat{\mu}^{S_n}$ denote the integer-valued random measure associated to the jumps of $\hat{X}^{S_n}$ (see [11], II.1.16). By Proposition 3.2, its $((\mathcal{F}_{t+}^X), \mathbb{P})$-compensator is $\nu^{S_n}$. It can easily be checked that

$$\frac{1}{3} \leq \frac{y \log y - y + 1}{(y-1)^2} \leq 1 \quad \text{for } y \in (0, 2]$$

and

$$\frac{1}{3} \leq \frac{y \log y - y + 1}{y - 1} \quad \text{for } y \geq 2.$$

(Notice however that $\lim_{y \to \infty} \frac{y \log y - y + 1}{y - 1} = \infty$.) Hence, it follows from (2.4) that

$$\mathbb{E}_\mathbb{P}[([\psi(X, \xi) - 1]^2 \wedge |\psi(X, \xi) - 1|) * \nu^{S_n}] < \infty$$

for the nonnegative measurable function $\psi: U \times \mathbb{R}^d \to \mathbb{R}_+$ defined by

$$\psi(x,\xi) := \phi_2(x)\mathbb{1}_{\{x+\xi=\partial\}} + \phi_3(x,\xi)\mathbb{1}_{\{x+\xi\in E\}}.$$

Consequently, by [11], II.1.33 c,

$$[(\psi(X_-,\xi) - 1)] * (\hat{\mu}^{S_n} - \nu^{S_n})$$

is a well defined $((\mathcal{F}^X_{t+}), \mathbb{P})$-local martingale. Moreover, it follows from (2.4) that

$$\mathbb{E}_{\mathbb{P}}\left[\int_0^{S_n} \langle \alpha(X_s)\phi_1(X_s), \phi_1(X_s)\rangle\, ds\right] < \infty.$$

Hence, by [11], III.4.5,

$$\phi_1(X) \cdot \hat{X}^{S_n,c}$$

is a well-defined continuous $((\mathcal{F}^X_{t+}), \mathbb{P})$-local martingale, where $\hat{X}^{S_n,c}$ denotes the continuous martingale part of $\hat{X}^{S_n}$, relative to the measure $\mathbb{P}$. In summary,

$$(4.1) \qquad L^n := \phi_1(X) \cdot \hat{X}^{S_n,c} + [(\psi(X_-,\xi) - 1)] * (\hat{\mu}^{S_n} - \nu^{S_n})$$

is a well-defined $((\mathcal{F}^X_{t+}), \mathbb{P})$-local martingale with

$$\langle L^{n,c}, L^{n,c}\rangle_\infty = \langle L^{n,c}, L^{n,c}\rangle_{S_n} = \int_0^{S_n} \langle \alpha(X_s)\phi_1(X_s), \phi_1(X_s)\rangle\, ds$$

and

$$\Delta L^n_t = [\psi(X_{t-}, \Delta\hat{X}_t) - 1]\mathbb{1}_{\{\Delta\hat{X}^{S_n}_t \neq 0\}} > -1.$$

This latter property assures that the stochastic exponential $\mathcal{E}(L^n)$ is a strictly positive $((\mathcal{F}^X_{t+}), \mathbb{P})$-local martingale. Moreover, it follows from Théorème IV.3 of [19], together with (2.4), that $\mathcal{E}(L^n)$ is a uniformly integrable $((\mathcal{F}^X_{t+}), \mathbb{P})$-martingale, which implies that

$$(4.2) \qquad D^n := D_0 \mathcal{E}(L^n)$$

is a nonnegative, uniformly integrable $((\mathcal{F}^X_{t+}), \mathbb{P})$-martingale.

Obviously, for $n \geq m$,

$$D^n_t = D^m_t \qquad \text{for all } t \leq S_m.$$

Therefore, for $t < S_\infty(\omega)$, and also for $t = S_\infty(\omega)$ if $S_\infty(\omega) = S_m(\omega)$ for some $m \geq 1$, one can define

$$D^\infty_t(\omega) := \lim_{n\to\infty} D^n_t(\omega) \geq 0.$$

Note that, for all $n \geq 1$, $\mathcal{E}(L^n)$ is strictly positive. Hence,

$$(4.3) \qquad D^\infty_t > 0 \qquad \text{for all } t \in [0, S_\infty) \text{ on } \{D_0 > 0\}.$$

Since $(D_{S_n}^\infty)_{n \geq 1} = (D_{S_n}^n)_{n \geq 1}$ is a nonnegative martingale, the limit

$$D_{S_\infty}^\infty := \lim_{n \to \infty} D_{S_n}^\infty \geq 0$$

exists $\mathbb{P}$-almost surely, and

$$D_t := D_t^\infty \mathbb{1}_{\{t < S_\infty\}} + D_{S_\infty}^\infty \mathbb{1}_{\{S_\infty \leq t\}}, \qquad t \in [0, \infty],$$

is a nonnegative càdlàg process. It follows from Fatou's lemma that, for all $t \geq 0$ and every $(\mathcal{F}_{t+}^X)$-stopping time $S$,

$$\mathbb{E}_\mathbb{P}[D_S | \mathcal{F}_{t+}^X] \leq \lim_{n \to \infty} \mathbb{E}_\mathbb{P}[D_{S \wedge S_n} | \mathcal{F}_{t+}^X] = \lim_{n \to \infty} D_{t \wedge S \wedge S_n} = D_{t \wedge S}.$$

In particular, $D$ is a supermartingale and

(4.4) $$\mathbb{E}_\mathbb{P}[D_T] \leq 1.$$

Now, let $f \in C_c^2(E)$ and set $f(\partial) = f(\Delta) = 0$. Then, it follows from (3.2) that $M^{f,S_n} = \mathcal{M}^{f,S_n}$. By Remark 2.3 part 2, $M^{f,S_n}$ is an $((\mathcal{F}_{t+}^X), \mathbb{P})$-martingale, and obviously, it has bounded jumps. Therefore, it follows from Lemma III.3.14 in [11] that $\langle M^{f,S_n}, L^n \rangle$ and $\langle M^{f,S_n}, D^{S_n} \rangle$ exist and $\langle M^{f,S_n}, D^{S_n} \rangle = D_-^{S_n} \cdot \langle M^{f,S_n}, L^n \rangle$. It can be seen from II.2.36, II.2.43 and the proof of II.2.42 in [11] that

(4.5) $$M^{f,S_n} = \mathcal{M}^{f,S_n} = \nabla f(X) \cdot \hat{X}^{S_n,c} + [f(X_- + \xi) - f(X_-)] * (\hat{\mu}^{S_n} - \nu^{S_n})$$

is the decomposition of $M^{f,S_n}$ into a continuous and a purely discontinuous $((\mathcal{F}_{t+}^X), \mathbb{P})$-local martingale part. Hence,

$$\langle M^{f,S_n}, L^n \rangle_t = \int_0^t \langle \nabla f(X_s), \alpha(X_s) \phi_1(X_s) \rangle \, ds$$
$$+ ([f(X + \xi) - f(X)][\psi(X, \xi) - 1]) * \nu_t^{S_n},$$

which shows that, for all $t \geq 0$,

$$\tilde{M}_t^{f,S_n} := f(X_t) - f(X_0) - \int_0^t \tilde{\mathcal{A}} f(X_s) \, ds$$
$$= M_t^{f,S_n} - \langle M^{f,S_n}, L^n \rangle_t = M_t^{f,S_n} - \int_0^t \frac{1}{D_{s-}^{S_n}} \, d\langle M^{f,S_n}, D^{S_n} \rangle_s.$$

Thus, it follows from Girsanov's theorem for local martingales in the form of [11], III.3.11, that $\tilde{M}^{f,S_n}$ is an $((\mathcal{F}_{t+}^X), D^{S_n} \cdot \mathbb{P})$-martingale. By the definition of $\mathbb{Q}$ and the optional sampling theorem, $\tilde{M}^{f,S_n}$ is also an $((\mathcal{F}_t^X), \mathbb{Q})$-martingale. By Remark 2.3, part 3, we can apply Theorem 4.6.1 of [8] (observe that for the proof of [8], Theorem 4.6.1, it is only needed that $S_n$ is an $(\mathcal{F}_t^X)$-stopping time, see also [8], Lemma 4.5.16) to conclude that

$$D_{S_n} \cdot \mathbb{P} = \mathbb{Q} \qquad \text{on } \mathcal{F}_{S_n}^X.$$

Now, let $A \in \mathcal{F}_T^X$. It can easily be checked that, for all $n \geq 1$,

$$A \cap \{T < S_n\} \in \mathcal{F}_{S_n \wedge T}^X.$$

Thus,

$$\mathbb{Q}[A \cap \{T < S_\infty\}] = \lim_{n \to \infty} \mathbb{Q}[A \cap \{T < S_n\}] = \lim_{n \to \infty} \mathbb{E}_\mathbb{P}[D_{S_n \wedge T} \mathbb{1}_{\{T < S_n\}} \mathbb{1}_A]$$
(4.6)
$$= \lim_{n \to \infty} \mathbb{E}_\mathbb{P}[D_T \mathbb{1}_{\{T < S_n\}} \mathbb{1}_A] = \mathbb{E}_\mathbb{P}[D_T \mathbb{1}_{\{T < S_\infty\}} \mathbb{1}_A],$$

where the first and the last equality follow from the monotone convergence theorem. This proves (2.5).

Equation (4.6) applied to $A = \Omega$ yields

$$\mathbb{Q}[T < S_\infty] = \mathbb{E}_\mathbb{P}[D_T \mathbb{1}_{\{T < S_\infty\}}].$$

Hence, if $\mathbb{Q}[T < S_\infty] = 1$, then (4.4) shows that

(4.7) $\quad\quad \mathbb{E}_\mathbb{P}[D_T] = 1 \quad$ and $D_T = 0$ on $\{T \geq S_\infty\}$ $\mathbb{P}$-a.s.,

which proves (2.6).

If, in addition, $\mathbb{Q}|_{\mathcal{F}_0^X} \sim \mathbb{P}|_{\mathcal{F}_0^X}$, then $D_0 > 0$ $\mathbb{P}$-a.s. and it follows from (4.3) that $D_T > 0$ on $\{T < S_\infty\}$ $\mathbb{P}$-a.s., which together with (4.7) implies that

(4.8) $\quad\quad\quad \{T \geq S_\infty\} = \{D_T = 0\} \quad\quad \mathbb{P}$-a.s.

Property (2.7) is now a consequence of (4.8) and (2.6).

If $\mathbb{Q}[T < R_\infty] = 1$, then $\mathbb{Q}[T \wedge T_n < S_\infty] = 1$, for all $n \geq 1$. Therefore, it follows from (2.6) that

$$\mathbb{Q}|_{\mathcal{F}_{T \wedge T_n}^X} = D_{T \wedge T_n} \cdot \mathbb{P}|_{\mathcal{F}_{T \wedge T_n}^X}.$$

Moreover, since $\lim_{n \to \infty} T_n = T_{\text{expl}} \geq S_\infty$, we have $\lim_{n \to \infty} D_{T \wedge T_n} = D_T$ $\mathbb{P}$-a.s. Hence, if $(D_{T \wedge T_n})_{n \geq 1}$ is uniformly integrable, then $D_{T \wedge T_n} \to D_T$ in $L^1(\mathbb{P})$. Therefore,

$$\mathbb{Q}|_{\mathcal{F}_{T \wedge T_n}^X} = D_T \cdot \mathbb{P}|_{\mathcal{F}_{T \wedge T_n}^X},$$

for all $n \geq 1$, which, by Proposition 2.1, implies (2.8), and the theorem is proved.

**5. Carré-du-champ operator.** Part (2.8) of Theorem 2.4 yields absolute continuity of $\mathbb{Q}|_{\mathcal{F}_T^X}$ with respect to $\mathbb{P}|_{\mathcal{F}_T^X}$, also on $\{T \geq T_{\text{expl}}\}$. In this section we consider a special choice of $\phi_1$, $\phi_2$ and $\phi_3$, which even provides equivalence beyond explosion. This is an extension of [27], Section VIII.3, and involves the *carré-du-champ operator* $\Gamma : C_c^2(E) \times C_c^2(E) \to B(E)$ defined by

$$\Gamma(f, g) := \mathcal{A}(fg) - f\mathcal{A}g - g\mathcal{A}f.$$

In contrast to above, we now first introduce a probability measure $\mathbb{Q}$ such that $\mathbb{Q} \sim \mathbb{P}$ on $\mathcal{F}^X_{t+}$ for all $t \geq 0$, and then find the appropriate generator $\tilde{\mathcal{A}}$ for which $\mathbb{Q}$ solves the martingale problem.

Fix $h \in C^2_c(E)$. Then $H := e^h - 1 \in C^2_c(E)$, and we can define

$$D_t := e^{h(X_t) - h(X_0)} \exp\left(-\int_0^t \frac{\mathcal{A}H(X_s)}{e^{h(X_s)}} \, ds\right).$$

Integration by parts, using $d(e^{h(X)}) = dM^H + \mathcal{A}H(X)\, dt$, yields

$$(5.1) \quad dD_t = e^{-h(X_0)} \exp\left(-\int_0^t \frac{\mathcal{A}H(X_s)}{e^{h(X_s)}} \, ds\right) dM^H_t = D_{t-} e^{-h(X_{t-})} \, dM^H_t.$$

Since $D$ is uniformly bounded on compact time intervals, we conclude that $D$ is a strictly positive $((\mathcal{F}^X_{t+}), \mathbb{P})$-martingale. As in [26], Theorem IV.38.9, it can be deduced from the Daniell–Kolmogorov extension theorem that there exists a probability measure $\mathbb{Q}$ on $\mathcal{F}^X$ such that $\mathbb{Q} = D_t \cdot \mathbb{P}$ on $\mathcal{F}^X_{t+}$ for all $t \geq 0$.

In view of (3.2) [we set $H(\partial) = 0$] and (4.5), we have

$$M^{H,S_n} = \mathcal{M}^{H,S_n} = \nabla e^{h(X)} \cdot \hat{X}^{S_n,c} + (e^{h(X_- + \xi)} - e^{h(X_-)}) * (\hat{\mu}^{S_n} - \nu^{S_n}),$$

so that, together with (5.1), we obtain

$$D^{S_n} = \mathcal{E}(\nabla h(X) \cdot \hat{X}^{S_n,c} + (e^{h(X_- + \xi) - h(X_-)} - 1) * (\hat{\mu}^{S_n} - \nu^{S_n})),$$

for all $n \geq 1$. Comparing this to (4.1) suggests that we are in the situation of Theorem 2.4 with

$$(5.2) \quad \phi_1(x) = \nabla h(x), \quad \phi_2(x) = e^{-h(x)} \quad \text{and} \quad \phi_3(x,\xi) = e^{h(x+\xi) - h(x)},$$

which clearly satisfy (2.9)–(2.11) for all $x \in E$ and a fixed constant $K > 0$.

THEOREM 5.1. $\mathbb{Q}$ *is a solution of the martingale problem for* $\tilde{\mathcal{A}} : C^2_c(E) \to B(E)$ *given by*

$$(5.3) \qquad \tilde{\mathcal{A}}f := \mathcal{A}f + \frac{\Gamma(H,f)}{e^h},$$

*which equals* (2.2) *with* (2.3) *and* (5.2).

PROOF. A straightforward calculation yields

$$\Gamma(f,g)(x) = \langle \alpha(x) \nabla f(x), \nabla g(x) \rangle + \gamma(x) f(x) g(x)$$
$$+ \int_{\mathbb{R}^d} (f(x+\xi) - f(x))(g(x+\xi) - g(x)) \mu(x, d\xi),$$

which makes it easy to see that (5.3) equals (2.2) with (2.3) and (5.2).

Let $f \in C_c^2(E)$. Lemma 5.2 below shows that

$$\langle M^f, M^H \rangle_t = \int_0^t \Gamma(f, H)(X_s) \, ds, \qquad t \geq 0.$$

Therefore,

$$\tilde{M}_t^f := f(X_t) - f(X_0) - \int_0^t \tilde{\mathcal{A}} f(X_s) \, ds$$

$$= f(X_t) - f(X_0) - \int_0^t \mathcal{A} f(X_s) \, ds - \int_0^t e^{-h(X_s)} \, d\langle M^f, M^H \rangle_s$$

$$= M_t^f - \int_0^t \frac{1}{D_{s-}} \, d\langle M^f, D \rangle_s,$$

and it follows from Girsanov's theorem for local martingales [11], III.3.11, that $\tilde{M}^f$ is an $((\mathcal{F}_{t+}^X), \mathbb{Q})$-martingale, which proves the theorem. □

LEMMA 5.2. *If $f, g \in C_c^2(E)$, then*

$$\langle M^f, M^g \rangle_t = \int_0^t \Gamma(f, g)(X_s) \, ds.$$

PROOF. Literally the same as the proof of Proposition VIII.3.3 in [27]. □

**6. Example.** We here apply Theorem 2.4 to a one-dimensional diffusion with compound Poisson jumps and a constant killing rate. In [2], it is applied to a multi-dimensional diffusion, and in [1] to a multi-dimensional jump-diffusion model.

Let $(\Omega', \mathcal{F}', \mathbb{P}')$ be a probability space that carries the following three independent random objects: a one-dimensional standard Brownian motion $(W_t)_{t \geq 0}$; a compound Poisson process $(N_t)_{t \geq 0}$ with jump arrival rate $\lambda > 0$ and positive jumps that are distributed according to a probability measure $m$ on $(0, \infty)$; and an exponentially distributed random variable $\tau$ with mean $\frac{1}{\gamma} > 0$. Let $b_0 \geq 0$, $b_1 \in \mathbb{R}$ and $\sigma > 0$. It is well known that the SDE

(6.1) $$dV_t = (b_0 + b_1 V_t) \, dt + \sigma \sqrt{V_t} \, dW_t, \qquad V_0 = v > 0,$$

has a unique strong solution, $V$ stays nonnegative, and

(6.2) $$V \text{ never reaches zero if } b_0 \geq \frac{\sigma^2}{2}.$$

(Cox, Ingersoll and Ross [3] model the short term interest rates by the solution of an SDE of the form (6.1).) It follows from a comparison argument that the same is true for the equation

(6.3) $$dY_t = (b_0 + b_1 Y_t) \, dt + \sigma \sqrt{Y_t} \, dW_t + dN_t, \qquad Y_0 = y > 0.$$

The process
$$Z := Y \mathbb{1}_{[0,\tau)} + \Delta \mathbb{1}_{[\tau,\infty)}$$
takes values in $E_\Delta$, for $E = \mathbb{R}_+$, and its distribution $\mathbb{P}$ is a probability measure on the measurable space $(\Omega, \mathcal{F}^X)$ introduced in Section 2. It can be checked that $\mathbb{P}$ is a solution of the martingale problem for
$$\mathcal{A}f(x) = \tfrac{1}{2}\sigma^2 x f''(x) + (b_0 + b_1 x) f'(x) - \gamma f(x)$$
$$+ \int_0^\infty [f(x+\xi) - f(x)] \lambda m(d\xi).$$

Let
$$\tilde{\mathcal{A}}f(x) = \tfrac{1}{2}\sigma^2 x f''(x) + (\tilde{b}_0 + \tilde{b}_1 x) f'(x) - \tilde{\gamma}(x) f(x)$$
$$+ \int_0^\infty [f(x+\xi) - f(x)] \tilde{\mu}(x, d\xi),$$
where $\tilde{b}_0 \geq \frac{\sigma^2}{2}$, $\tilde{b}_1 \in \mathbb{R}$, $\tilde{\gamma}(x) = \tilde{\gamma}_0 + \tilde{\gamma}_1 x$, for some $(\tilde{\gamma}_0, \tilde{\gamma}_1) \in \mathbb{R}_+^2 \setminus \{(0,0)\}$, and $\tilde{\mu}(x, \cdot)$ is, for all $x > 0$, a measure on $(0, \infty)$ of the form $\tilde{\mu}(x, d\xi) = [m_0(\xi) + m_1(\xi)x]\lambda m(d\xi)$, for nonnegative measurable functions $m_0, m_1 : (0, \infty) \to \mathbb{R}_+$, such that $(m_0(\xi), m_1(\xi)) \in \mathbb{R}_+^2 \setminus \{(0,0)\}$ for all $\xi > 0$ and
$$\int_0^\infty l(m_0(\xi) + m_1(\xi)x) m(d\xi) < \infty \qquad \text{for all } x > 0,$$
where $l(u) = u \log u - u + 1$. It follows from Theorem 2.7 in [7] that the martingale problem for $\tilde{\mathcal{A}}$ is well posed. Let $\mathbb{Q}$ be the solution of the martingale problem for $\tilde{\mathcal{A}}$ with initial distribution $\delta_y$. It can be deduced from (6.2) and a comparison argument that

(6.4) $\quad \mathbb{Q}[\text{there exists a } t \geq 0 \text{ such that } X_t = 0 \text{ or } X_{t-} = 0] = 0.$

We set $U = (0, \infty)$ and $U_n = (1/n, n)$, $n \geq 1$. Since we have no explosion, (6.4) implies that $\mathbb{Q}[S_\infty = \infty] = 1$. Furthermore, the measurable mappings
$$\phi_1(x) = \frac{\tilde{b}_0 - b_0}{\sigma^2 x} + \frac{\tilde{b}_1 - b_1}{\sigma^2}, \qquad x \in U,$$
$$\phi_2(x) = \frac{1}{\gamma}(\tilde{\gamma}_0 + \tilde{\gamma}_1 x), \qquad x \in U,$$
$$\phi_3(x, \xi) = \begin{cases} m_0(\xi) + m_1(\xi)x, & \text{if } \xi > 0, \\ 1, & \text{if } \xi \leq 0, \end{cases} x \in U,$$
satisfy the conditions (2.9)–(2.11), and for all $x \in U$,
$$\tilde{b}_0 + \tilde{b}_1 x = b_0 + b_1 x + \sigma^2 x \phi_1(x),$$
$$\tilde{\gamma}(x) = \phi_2(x)\gamma,$$
$$\tilde{\mu}(x, d\xi) = \phi_3(x, \xi)\lambda m(d\xi).$$

Therefore, Theorem 2.4 applies, and we obtain that

$$\mathbb{Q}|_{\mathcal{F}_T^X} \ll \mathbb{P}|_{\mathcal{F}_T^X}$$

for all $(\mathcal{F}_t^X)$-stopping times $T < \infty$. Moreover, if $b_0 \geq \frac{\sigma^2}{2}$, then

$$\mathbb{P}[S_\infty = \infty] = 1 - \mathbb{P}[\text{there exists a } t \geq 0 \text{ such that } X_t = 0 \text{ or } X_{t-} = 0] = 1,$$

and Theorem 2.4 yields that

$$\mathbb{Q}|_{\mathcal{F}_T^X} \sim \mathbb{P}|_{\mathcal{F}_T^X}$$

for all $(\mathcal{F}_t^X)$-stopping times $T < \infty$. If we identify $\Delta$ with $-1$, the process $Z$ becomes the semimartingale

$$\hat{Z} := Y \mathbb{1}_{[0,\tau)} - \mathbb{1}_{[\tau,\infty)}.$$

It can be seen from (6.3) that $d\hat{Z}_t^c = \mathbb{1}_{\{0 \leq t < \tau\}} \sigma \sqrt{Y_t}\, dW_t$. The random measure $\hat{\mu}$ associated to the jumps of $\hat{Z}$ is an integer-valued random measure on $\mathbb{R}_+^2$ with compensator

$$\nu(dt, d\xi) = \mathbb{1}_{\{0 \leq t < \tau\}}\, dt \times (\lambda m(d\xi) + \gamma \delta_{-1-Z_t}(d\xi)).$$

Since the distribution of

$$\psi(Z_{t-}, \xi) = \phi_2(Z_{t-}) \mathbb{1}_{\{Z_{t-}+\xi=-1\}} + \phi_3(Z_{t-}, \xi) \mathbb{1}_{\{Z_{t-}+\xi \geq 0\}},$$

and the stochastic exponential

$$D' = \mathcal{E}(\phi_1(Z) \cdot \hat{Z}^c + (\psi(Z_-, \xi) - 1) * (\hat{\mu} - \nu))$$

only depend on the distribution of $Z$, it follows from (4.7) that

$$\mathbb{E}_{\mathbb{P}'}[D_t'] = \mathbb{E}_\mathbb{P}[D_t] = 1.$$

Hence, $D'$ is a $\mathbb{P}'$-martingale, and for all $t \geq 0$, $D_t' \cdot \mathbb{P}'$ is a probability measure on $(\Omega', \mathcal{F}')$ under which the distribution of the stopped process $Z^t$ is equal to $\mathbb{Q}|_{\mathcal{F}_t^X}$. If $b_0 \geq \frac{\sigma^2}{2}$, then $D_t' > 0$ $\mathbb{P}'$-almost surely for all $t \in \mathbb{R}_+$, and $D_t' \cdot \mathbb{P}'$ is equivalent to $\mathbb{P}'$.

## APPENDIX

PROOF OF PROPOSITION 2.1. It is clear that

(A.1) $$\mathcal{F}_T^X \supset \mathcal{F}_{T \wedge T_{\text{expl}}}^X \supset \sigma\left(\bigcup_{n \geq 1} \mathcal{F}_{T \wedge T_n}^X\right).$$

To show the reverse inclusions, we first prove that

(A.2) $$\mathcal{F}^X \subset \sigma\left(\bigcup_{n \geq 1} \mathcal{F}_{T_n}^X\right).$$

Note that for all $t \geq 0$, and all Borel subsets $B$ of $E$,

$$\{X_t \in B\} = \{X_t \in B\} \cap \{T_{\text{expl}} > t\} = \bigcup_{n \geq 1} (\{X_t \in B\} \cap \{T_n > t\}),$$

and for all $n \geq 1$,

$$\{X_t \in B\} \cap \{T_n > t\} \in \mathcal{F}^X_{T_n}.$$

Hence,

(A.3) $$\{X_t \in B\} \in \sigma\left(\bigcup_{n \geq 1} \mathcal{F}^X_{T_n}\right).$$

Moreover, for all $t \geq 0$,

$$\{X_t = \Delta\} = \{T_{\text{expl}} \leq t\} \cup \{T_{\text{jump}} \leq t\}$$

$$= \left(\bigcap_{n \geq 1} \{T_n \leq t\}\right) \cup (\{T_{\text{jump}} \leq t\} \cap \{T_{\text{expl}} > t\})$$

$$= \left(\bigcap_{n \geq 1} \{T_n \leq t\}\right) \cup \bigcup_{n \geq 1} (\{T_{\text{jump}} \leq t\} \cap \{T_n > t\}).$$

It can easily be checked that, for all $n \geq 1$,

$$\{T_n \leq t\} \quad \text{and} \quad \{T_{\text{jump}} \leq t\} \cap \{T_n > t\} \text{ belong to } \mathcal{F}^X_{T_n}.$$

Hence,

$$\{X_t = \Delta\} \in \sigma\left(\bigcup_{n \geq 1} \mathcal{F}^X_{T_n}\right),$$

which, together with (A.3), implies (A.2).

For every set $A \in \mathcal{F}^X_T$, we write

(A.4) $$A = [A \cap \{T < T_{\text{expl}}\}] \cup [A \cap \{T \geq T_{\text{expl}}\}]$$

$$= \left[\bigcup_{n \geq 1} A \cap \{T < T_n\}\right] \cup [A \cap \{T \geq T_{\text{expl}}\}].$$

Observe that, for all $n \geq 1$,

(A.5) $$A \cap \{T < T_n\} \in \mathcal{F}^X_{T \wedge T_n}.$$

For every class of subset $\mathcal{G}$ of $\Omega$, we define

$$\mathcal{G} \cap \{T \geq T_{\text{expl}}\} := \{G \cap \{T \geq T_{\text{expl}}\} \mid G \in \mathcal{G}\}.$$

It follows from (A.2) that

$$(A.6) \qquad A \cap \{T \geq T_{\text{expl}}\} \in \sigma\left(\bigcup_{n \geq 1} \mathcal{F}^X_{T_n}\right) \cap \{T \geq T_{\text{expl}}\},$$

and it can easily be checked that

$$\sigma\left(\bigcup_{n \geq 1} \mathcal{F}^X_{T_n}\right) \cap \{T \geq T_{\text{expl}}\} \subset \sigma\left(\bigcup_{n \geq 1} \mathcal{F}^X_{T_n} \cap \{T \geq T_{\text{expl}}\}\right) \subset \sigma\left(\bigcup_{n \geq 1} \mathcal{F}^X_{T \wedge T_n}\right).$$

Hence, (A.4), (A.5) and (A.6) imply that

$$\mathcal{F}^X_T \subset \sigma\left(\bigcup_{n \geq 1} \mathcal{F}^X_{T \wedge T_n}\right),$$

which, together with (A.1), proves the proposition. $\square$

**Acknowledgments.** We are grateful to the referees and the Associate Editor for many valuable comments.

P. CHERIDITO
DEPARTMENT OF OPERATIONS RESEARCH
  AND FINANCIAL ENGINEERING
PRINCETON UNIVERSITY
PRINCETON, NEW JERSEY 08544
USA
E-MAIL: dito@princeton.edu

D. FILIPOVIĆ
DEPARTMENT OF MATHEMATICS
UNIVERSITY OF MUNICH
THERESIENSTRASSE 39
80333 MUNICH
GERMANY
E-MAIL: filipo@math.lmu.de

M. Yor
Laboratoire de Probabilités
  et Modèles aléatoires
Université Pierre et Marie Curie, Paris VI
4 Place Jussieu
75252 Paris Cedex 05
France